\documentclass[12pt]{article}%
\usepackage{amsmath}
\usepackage{amsfonts}
\usepackage{amssymb}
\usepackage{graphicx}%
\newtheorem{theorem}{Theorem}
\newtheorem{acknowledgement}[theorem]{Acknowledgement}

\newtheorem{remark}[theorem]{Remark}

\newenvironment{proof}[1][Proof]{\noindent\textbf{#1.} }{\ \rule{0.5em}{0.5em}}
\begin{document}

\title{On a complex differential Riccati equation}
\author{Kira V. Khmelnytskaya and Vladislav V. Kravchenko\\{\small Department of Mathematics, CINVESTAV del IPN, Unidad Quer\'{e}taro}\\{\small Libramiento Norponiente No. 2000, Fracc. Real de Juriquilla}\\{\small Queretaro, Qro. C.P. 76230 MEXICO}\\{\small e-mail: vkravchenko@qro.cinvestav.mx}}
\maketitle

\begin{abstract}
We consider a nonlinear partial differential equation for complex-valued
functions which is related to the two-dimensional stationary Schr\"{o}dinger
equation and enjoys many properties similar to those of the ordinary
differential Riccati equation such as the famous Euler theorems, the Picard
theorem and others. Besides these generalizations of the classical
\textquotedblleft one-dimensional\textquotedblright\ results we discuss new
features of the considered equation including an analogue of the Cauchy
integral theorem.

\end{abstract}

\section{Introduction}

The ordinary differential Riccati equation%
\index{Riccati ordinary differential equation}
\begin{equation}
u^{\prime}=pu^{2}+qu+r \label{Riccati}%
\end{equation}
has received a great deal of attention since a particular version was first
studied by Count Riccati in 1724, owing to both its specific properties and
the wide range of applications in which it appears. For a survey of the
history and classical results on this equation, see for example \cite{davis}
and \cite{reid}. \ This equation can be always reduced to its canonical form
(see, e.g., \cite{bogd}, \cite{Kamke}),
\begin{equation}
y^{\prime}+y^{2}=v, \label{ricc1}%
\end{equation}
and this is the form that we will consider.

One of the reasons for which the Riccati equation has so many applications is
that it is related to the general second order homogeneous differential
equation. \ In particular, the one-dimensional Schr\"{o}dinger equation
\begin{equation}
-\partial^{2}u+vu=0 \label{schrod1}%
\end{equation}
is related to (\ref{ricc1}) by the easily inverted substitution
\[
y=\frac{u^{\prime}}{u}.
\]
\ This substitution, which as its most spectacular application reduces
Burger's equation to the standard one-dimensional heat equation, is the basis
of the well-developed theory of logarithmic derivatives for the integration of
nonlinear differential equations \cite{march}. \ A generalization of this
substitution will be used in this work.

A second relation between the one-dimensional Schr\"{o}dinger equation and the
Riccati equation is as follows. \ The one-dimensional Schr\"{o}dinger operator
can be factorized in the form
\begin{equation}
-\frac{d^{2}}{dx^{2}}+v(x)=-(\frac{d}{dx}+y(x))(\frac{d}{dx}-y(x))
\label{fact1}%
\end{equation}
if and only if (\ref{ricc1}) holds. This observation is a key to a vast area
of research related to the factorization method (see, e.g.,
\cite{MielnikReyes} and \cite{MielnikRosas2004}) and to the Darboux
transformation (see, e.g., \cite{MS}, \cite{ND} and \cite{Samsonov95}). In the
present work we consider a result similar to (\ref{fact1}) but already in a
two-dimensional setting.

Among the peculiar properties of the Riccati equation stand out two theorems
of Euler, dating from 1760. \ The first of these states that if a particular
solution $y_{0}$ of the Riccati equation is known, the substitution
$y=y_{0}+z$ reduces (\ref{ricc1}) to a Bernoulli equation which in turn is
reduced by the substitution $z=\frac{1}{u}$ to a first order linear equation.
\ Thus given a particular solution of the Riccati equation, it can be
linearized and the general solution can be found in two integrations. \ The
second of these theorems states that given two particular solutions
$y_{0},\ y_{1}$ of the Riccati equation, the general solution can be found in
the form
\begin{equation}
y=\frac{ky_{0}\exp(\int y_{0}-y_{1})-y_{1}}{k\exp(\int y_{0}-y_{1})-1}
\label{2sol}%
\end{equation}
where $k$ is a constant. \ That is, given two particular solutions of
(\ref{ricc1}), the general solution can be found in one integration.

Other interesting properties are those discovered by Weyr and Picard
(\cite{wats}, \cite{davis}). \ The first is that given a third particular
solution $y_{3}$, the general solution can be found without integrating.
\ That is, an explicit combination of three particular solutions gives the
general solution. The second is that given a fourth particular solution
$y_{4}$, the cross ratio
\[
\frac{(y_{1}-y_{2})(y_{3}-y_{4})}{(y_{1}-y_{4})(y_{3}-y_{2})}%
\]
is a constant. Thus the derivative of this ratio is zero, which holds if and
only if the numerator of the derivative is zero, that is, if and only if%
\[
(y_{1}-y_{4})(y_{3}-y_{2})((y_{1}-y_{2})(y_{3}-y_{4}))^{\prime}-(y_{1}%
-y_{2})(y_{3}-y_{4})((y_{1}-y_{4})(y_{3}-y_{2}))^{\prime}=0.
\]

Dividing by $(y_{1}-y_{2})(y_{3}-y_{4})(y_{1}-y_{4})(y_{3}-y_{2})$, we see
that Picard's theorem is equivalent to the statement
\begin{equation}
\frac{(y_{1}-y_{2})^{\prime}}{y_{1}-y_{2}}+\frac{(y_{3}-y_{4})^{\prime}}%
{y_{3}-y_{4}}-\frac{(y_{1}-y_{4})^{\prime}}{y_{1}-y_{4}}-\frac{(y_{3}%
-y_{2})^{\prime}}{y_{3}-y_{2}}=0. \label{weyrgen}%
\end{equation}

In the present work we study the following equation
\begin{equation}
\partial_{\overline{z}}Q+\left\vert Q\right\vert ^{2}=v \label{ComplexRiccati}%
\end{equation}
where $z$ is a complex variable, $\partial_{\overline{z}}=\frac{1}{2}%
(\partial_{x}+i\partial_{y})$, $Q$ is a complex valued function of $z$ and $v$
is a real valued function. Note that this equation is different from the
complex Riccati equation studied in dozens of works where it is supposed to
have the form (\ref{Riccati}) with complex analytic coefficients $p$, $q$ and
$r$ and a complex analytic solution $u$ (see, e.g., \cite{Hille}). We do not
suppose analyticity of the complex functions involved and show that equation
(\ref{ComplexRiccati}) unlike the equation considered in \cite{Hille} is
related to the two-dimensional stationary Schr\"{o}dinger equation in a
similar way as the ordinary Riccati and Schr\"{o}dinger equations are related
in dimension one. Moreover, we establish generalizations of the Euler and
Picard theorems and obtain some other results which are essentially
two-dimensional, for example, an analogue of the Cauchy integral theorem for
solutions of the complex Riccati equation (\ref{ComplexRiccati}).

Equation (\ref{ComplexRiccati}) first appeared in \cite{KrBers} as a reduction
to a two-dimensional case of the spatial factorization of the stationary
Schr\"{o}dinger operator which was studied in a quaternionic setting in
\cite{Swansolo}, \cite{Swan}, \cite{KK}, \cite{KKW}, \cite{AQA}, \cite{KrBers}
and later on using Clifford analysis in \cite{Bernstein2006} and
\cite{DeSchepperPena}.

The ordinary Riccati equation is at the heart of many analytical and numerical
approaches to problems involving the one-dimensional Schr\"{o}dinger and
Sturm-Liouville equations. Here we furnish a complete structural analogy
between dimensions one and two regarding the relationship between the
Schr\"{o}dinger and the Riccati equations. Besides, the deep similarity
between the ordinary Riccati equation and (\ref{ComplexRiccati}) strongly
suggests that many known applications of the ordinary Riccati equation can be
generalized to the two-dimensional situation and many new aspects such as
theorem \ref{Th_Cauchy_Int_theorem_Riccati} will become manifest.

\section{Some preliminary results on the stationary Schr\"{o}dinger equation
and a class of pseudoanalytic functions}

We need the following definition. Consider the equation $\partial
_{\overline{z}}\varphi=\Phi$ in the whole complex plane or in a convex domain,
where $\Phi=\Phi_{1}+i\Phi_{2}$ is a given complex valued function whose real
part $\Phi_{1}$ and imaginary part $\Phi_{2}$ satisfy the equation
\begin{equation}
\partial_{y}\Phi_{1}-\partial_{x}\Phi_{2}=0, \label{casirot}%
\end{equation}
then we can reconstruct $\varphi$ up to an arbitrary real constant $c$ in the
following way%
\begin{equation}
\varphi(x,y)=2\left(  \int_{x_{0}}^{x}\Phi_{1}(\eta,y)d\eta+\int_{y_{0}}%
^{y}\Phi_{2}(x_{0},\xi)d\xi\right)  +c \label{Antigr}%
\end{equation}
where $(x_{0},y_{0})$ is an arbitrary fixed point in the domain of interest.

By $\overline{A}$ we denote the integral operator in (\ref{Antigr}):%
\[
\overline{A}[\Phi](x,y)=2\left(  \int_{x_{0}}^{x}\Phi_{1}(\eta,y)d\eta
+\int_{y_{0}}^{y}\Phi_{2}(x_{0},\xi)d\xi\right)  +c.
\]
Note that formula (\ref{Antigr}) can be easily extended to any simply
connected domain by considering the integral along an arbitrary rectifiable
curve $\Gamma$ leading from $(x_{0},y_{0})$ to $(x,y),$%
\[
\varphi(x,y)=2\left(  \int_{\Gamma}\Phi_{1}dx+\Phi_{2}dy\right)  +c.
\]
Thus if $\Phi$ satisfies (\ref{casirot}), there exists a family of real valued
functions $\varphi$ such that $\partial_{\overline{z}}\varphi=\Phi$, given by
the formula $\varphi=\overline{A}[\Phi]$. In a similar way we introduce the
operator
\[
A[\Phi](x,y)=2\left(  \int_{\Gamma}\Phi_{1}dx-\Phi_{2}dy\right)  +c
\]
which is applicable to complex functions satisfying the condition
\begin{equation}
\partial_{y}\Phi_{1}+\partial_{x}\Phi_{2}=0 \label{casirot+}%
\end{equation}
and corresponds to the operator $\partial_{z}$.

Let $f$ denote a positive twice continuously differentiable function defined
in a domain $\Omega\subset\mathbb{C}$. Consider the following Vekua equation%
\begin{equation}
W_{\overline{z}}=\frac{f_{\overline{z}}}{f}\overline{W}\qquad\text{in }%
\Omega\label{Vekuamain}%
\end{equation}
where the subindex $\overline{z}$ means the application of the operator
$\partial_{\overline{z}}$ (analogously, with the aid of the subindex $z$ we
will denote the application of the operator $\partial_{z}$), $W$ is a
complex-valued function and $\overline{W}=C[W]$ is its complex conjugate
function. As was shown in \cite{KrBers}, \cite{Krpseudoan} (see also
\cite{KrJPhys06}, \cite{KrRecentDevelopments}, \cite{KrOviedo06}) equation
(\ref{Vekuamain}) is closely related to the second-order equation of the form
\begin{equation}
\left(  -\Delta+\nu\right)  u=0\qquad\text{in }\Omega\label{SchrIntr}%
\end{equation}
where $\nu=\left(  \Delta f\right)  /f$ and $u$ are real-valued functions. In
particular the following statements are valid.

\begin{theorem}
\cite{Krpseudoan} \label{PrDarboux}Let $W$ be a solution of (\ref{Vekuamain}).
Then $u=\operatorname{Re}W$ is a solution of (\ref{SchrIntr}) and
$v=\operatorname{Im}W$ is a solution of the equation%
\begin{equation}
\left(  -\Delta+\eta\right)  v=0 \label{SchrodDarboux}%
\end{equation}
where $\eta=2\left(  \frac{\left\vert \nabla f\right\vert }{f}\right)
^{2}-\nu$.
\end{theorem}

\begin{theorem}
\cite{Krpseudoan} \label{PrTransform}Let $u$ be a solution of (\ref{SchrIntr})
in a simply connected domain $\Omega$. Then the function
\[
v\in\ker\left(  \Delta+\nu-2\left(  \frac{\left\vert \nabla f\right\vert }%
{f}\right)  ^{2}\right)
\]
such that $W=u+iv$ be a solution of (\ref{Vekuamain}), is constructed
according to the formula%
\begin{equation}
v=f^{-1}\overline{A}(if^{2}\partial_{\overline{z}}(f^{-1}u)).
\label{transfDarboux}%
\end{equation}
It is unique up to an additive term $cf^{-1}$ where $c$ is an arbitrary real constant.

Given $v\in\ker\left(  \Delta+\nu-2\left(  \frac{\left\vert \nabla
f\right\vert }{f}\right)  ^{2}\right)  ,$ the corresponding $u\in\ker\left(
\Delta-\nu\right)  $ can be constructed as follows%
\begin{equation}
u=-f\overline{A}(if^{-2}\partial_{\overline{z}}(fv)) \label{transfDarbouxinv}%
\end{equation}
up to an additive term $cf.$
\end{theorem}

Thus, the relation between (\ref{Vekuamain}) and (\ref{SchrIntr}) is very
similar to that between the Cauchy-Riemann system and the Laplace equation.
Moreover, choosing $\nu\equiv0$ and $f\equiv1$ we arrive at the well known
formulas from classical complex analysis. Note that the potential $\eta$ has
the form of a potential obtained after the Darboux transformation but already
in two dimensions.

For a Vekua equation of the form%
\[
W_{\overline{z}}=aW+b\overline{W}%
\]
where $a$ and $\ b$ are arbitrary complex-valued functions from an appropriate
function space \cite{Vekua}, a well developed theory of Taylor and Laurent
series in formal powers has been created (see \cite{Berskniga},
\cite{BersFormalPowers}). We recall that a formal power $Z^{(n)}(a,z_{0};z)$
corresponding to a coefficient $a$ and a centre $z_{0}$ is a solution of the
Vekua equation satisfying the asymptotic formula
\begin{equation}
Z^{(n)}(a,z_{0};z)\sim a(z-z_{0})^{n},\quad z\rightarrow z_{0}.
\label{asymptformulas}%
\end{equation}
For a rigorous definition we refer to \cite{Berskniga},
\cite{BersFormalPowers}. The theory of Taylor and Laurent series in formal
powers among other results contains the expansion and the Runge theorems as
well as more precise convergence results \cite{Menke} and a recently obtained
simple algorithm \cite{KrRecentDevelopments} for explicit construction of
formal powers for the Vekua equation of the form (\ref{Vekuamain}). In section
\ref{SectGeneralizations} we show how this theory is applied for generalizing
the second Euler theorem. For this we need the expansion theorem from
\cite{Berskniga}. For the Vekua equation of the form (\ref{Vekuamain}) this
expansion theorem reads as follows (for more details we refer the reader to
\cite{KrRecentDevelopments}).

\begin{theorem}
\label{ThTaylorRepr}Let $W$ be a regular solution of (\ref{Vekuamain}) defined
for $\left\vert z-z_{0}\right\vert <R$. Then it admits a unique expansion of
the form $W(z)=\sum_{n=0}^{\infty}Z^{(n)}(a_{n},z_{0};z)$ which converges
normally for $\left\vert z-z_{0}\right\vert <R$.
\end{theorem}

Another result which will be used in the present work (section
\ref{ThCauchySchr}) is a Cauchy-type integral theorem for the stationary
Schr\"{o}dinger equation. It was obtained in \cite{KrBers} with the aid of the
pseudoanalytic function theory.

\begin{theorem}
\label{ThCauchySchr}(Cauchy's integral theorem for the Schr\"{o}dinger
equation) Let $f$ be a nonvanishing solution of (\ref{SchrIntr}) in a domain
$\Omega$ and $u$ be another arbitrary solution of (\ref{SchrIntr}) in $\Omega
$. Then for every closed curve $\Gamma$ situated in a simply connected
subdomain of $\Omega$,%
\begin{equation}
\operatorname*{Re}\int_{\Gamma}\partial_{z}(\frac{u}{f})dz+i\operatorname*{Im}%
\int_{\Gamma}f^{2}\partial_{z}(\frac{u}{f})dz=0. \label{CauchySchrTh}%
\end{equation}

\end{theorem}

\section{The two-dimensional stationary Schr\"{o}dinger equation and the
complex Riccati equation}

Consider the complex differential Riccati equation
\begin{equation}
\partial_{\overline{z}}Q+\left\vert Q\right\vert ^{2}=\frac{\nu}{4}
\label{RiccatiMain}%
\end{equation}
where for convenience the factor $1/4$ was included. We recall that $\nu$ is a
real-valued function. Together with this equation we consider the stationary
Schr\"{o}dinger equation (\ref{SchrIntr}),
\begin{equation}
\left(  -\Delta+\nu\right)  u=0 \label{Schrod2}%
\end{equation}
where $u$ is real-valued. Both equations are studied in a domain
$\Omega\subset\mathbb{C}$.

\begin{theorem}
\label{ThFromSchrodToRiccati}Let $u$ be a solution of (\ref{Schrod2}). Then
its logarithmic derivative%
\begin{equation}
Q=\frac{u_{z}}{u} \label{logder}%
\end{equation}
is a solution of (\ref{RiccatiMain}).
\end{theorem}

\begin{proof}
It is only necessary to substitute (\ref{logder}) into (\ref{RiccatiMain}).
\end{proof}

\begin{remark}
Any solution of equation (\ref{RiccatiMain}) fulfils (\ref{casirot+}). Indeed,
the imaginary part of (\ref{RiccatiMain}) reads as follows%
\[
\partial_{y}Q_{1}+\partial_{x}Q_{2}=0.
\]

\end{remark}

\begin{theorem}
\label{ThFromRiccatiToSchrod}Let $Q$ be a solution of (\ref{RiccatiMain}).
Then the function
\begin{equation}
u=e^{A[Q]} \label{uexp}%
\end{equation}
is a solution of (\ref{Schrod2}).
\end{theorem}

\begin{proof}
Equation (\ref{Schrod2}) can be written in the form
\[
(4\partial_{\overline{z}}\partial_{z}-\nu)u=0.
\]
Taking $u$ in the form (\ref{uexp}) where $Q$ is a solution of
(\ref{RiccatiMain}) and using the observation that
\[
\partial_{\overline{z}}(A[Q])=\overline{\partial_{z}(A[Q])}=\overline{Q}%
\]
we have
\[
\partial_{\overline{z}}\partial_{z}u=\partial_{\overline{z}}\left(
Qe^{A[Q]}\right)  =e^{A[Q]}\left(  \partial_{\overline{z}}Q+\left\vert
Q\right\vert ^{2}\right)  =\frac{\nu}{4}u.
\]

\end{proof}

Observe that this theorem means that if $Q$ is a solution of
(\ref{RiccatiMain}) then there exists a solution $u$ of (\ref{Schrod2}) such
that (\ref{logder}) is valid. Theorems \ref{ThFromSchrodToRiccati} and
\ref{ThFromRiccatiToSchrod} are direct generalizations of the corresponding
facts from the one-dimensional theory.

The following statement is a generalization of the one-dimensional
factorization (\ref{fact1}).

\begin{theorem}
Given a complex function $Q$, for any real valued twice continuously
differentiable function $\varphi$ the following equality is valid
\begin{align}
\frac{1}{4}\left(  \Delta-\nu\right)  \varphi &  =(\partial_{\overline{z}%
}+QC)(\partial_{z}-QC)\varphi\label{fact2}\\
&  =(\partial_{z}+\overline{Q}C)(\partial_{\overline{z}}-\overline{Q}%
C)\varphi\nonumber
\end{align}
if and only if $Q$ is a solution of the Riccati equation (\ref{RiccatiMain}).
\end{theorem}

\begin{proof}
Consider%
\[
(\partial_{\overline{z}}+QC)(\partial_{z}-QC)\varphi=\frac{1}{4}\Delta
\varphi-\left\vert Q\right\vert ^{2}\varphi-Q_{\overline{z}}\varphi
\]
from which it is seen that (\ref{fact2}) is valid iff $Q$ is a solution of
(\ref{RiccatiMain}). The second equality in (\ref{fact2}) is obtained by
applying $C$ to both sides of the first equality.
\end{proof}

\section{Generalizations of classical theorems\label{SectGeneralizations}}

In this section we give generalizations of both Euler's theorems for the
Riccati equation as well as of Picard's theorem.

\begin{theorem}
(First Euler's theorem) Let $Q_{0}$ be a bounded particular solution of
(\ref{RiccatiMain}). Then (\ref{RiccatiMain}) reduces to the following first
order (real-linear) equation%
\begin{equation}
W_{\overline{z}}=\overline{Q_{0}W} \label{VekuaQ}%
\end{equation}
in the following sense. Any solution of (\ref{RiccatiMain}) has the form
\[
Q=\frac{\partial_{z}\operatorname*{Re}W}{\operatorname*{Re}W}%
\]
and vice versa, any solution of (\ref{VekuaQ}) can be expressed via a
corresponding solution $Q$ of (\ref{RiccatiMain}) as follows%
\begin{equation}
W=e^{A[Q]}+ie^{-A[Q_{0}]}\overline{A}\left[  ie^{2A[Q_{0}]}\partial
_{\overline{z}}e^{A[Q-Q_{0}]}\right]  . \label{WviaQ}%
\end{equation}

\end{theorem}

\begin{proof}
Let $Q_{0}$ be a bounded solution of (\ref{RiccatiMain}). Then by theorem
\ref{ThFromRiccatiToSchrod} we have that there exists a nonvanishing real
valued solution $f$ of (\ref{Schrod2}) such that $Q_{0}=f_{z}/f$. Hence
(\ref{VekuaQ}) has the form (\ref{Vekuamain}). Now, let $Q$ be any solution of
(\ref{RiccatiMain}). Then again $Q=u_{z}/u$ where $u$ is a solution of
(\ref{Schrod2}). According to theorem \ref{PrTransform}, $u$ is a real part of
a solution $W$ of (\ref{Vekuamain}). Thus we have proved the first part of the theorem.

Let $W=u+iv$ be any solution of (\ref{VekuaQ}) ($u=\operatorname*{Re}W$). Then
$u$ is a solution of (\ref{Schrod2}) and by theorem
\ref{ThFromSchrodToRiccati} it can be represented in the form $u=e^{A[Q]}$
where $Q$ is a solution of (\ref{RiccatiMain}). Then by theorem
\ref{PrTransform} (formula (\ref{transfDarboux})), $W$ has the form
(\ref{WviaQ}).
\end{proof}

Thus the Riccati equation (\ref{RiccatiMain}) is equivalent to a main Vekua
equation of the form (\ref{Vekuamain}).

In what follows we suppose that in the domain of interest $\Omega$ there
exists a bounded solution of (\ref{RiccatiMain}).

\begin{theorem}
(Second Euler's theorem) Any solution $Q$ of equation (\ref{RiccatiMain})
defined for $\left\vert z-z_{0}\right\vert <R$ can be represented in the form
\begin{equation}
Q=\frac{\sum_{n=0}^{\infty}Q_{n}e^{A[Q_{n}]}}{\sum_{n=0}^{\infty}e^{A[Q_{n}]}}
\label{QviaQ_n}%
\end{equation}
where $\left\{  Q_{n}\right\}  _{n=0}^{\infty}$ is the set of particular
solutions of the Riccati equation (\ref{RiccatiMain}) obtained as follows
\[
Q_{n}(z)=\frac{\partial_{z}\operatorname*{Re}Z^{(n)}(a_{n},z_{0}%
,z)}{\operatorname*{Re}Z^{(n)}(a_{n},z_{0},z)},
\]
$Z^{(n)}(a_{n},z_{0},z)$ are formal powers corresponding to equation
(\ref{VekuaQ}) and both series in (\ref{QviaQ_n}) converge normally for
$\left\vert z-z_{0}\right\vert <R$.
\end{theorem}

\begin{proof}
By the first Euler theorem we have
\[
Q=\frac{\partial_{z}\operatorname*{Re}W}{\operatorname*{Re}W}%
\]
where $W$ is a solution of (\ref{VekuaQ}). From theorem \ref{ThTaylorRepr} we
obtain%
\[
Q(z)=\frac{\partial_{z}\sum_{n=0}^{\infty}\operatorname*{Re}Z^{(n)}%
(a_{n},z_{0};z)}{\sum_{n=0}^{\infty}\operatorname*{Re}Z^{(n)}(a_{n},z_{0}%
;z)}.
\]
Every formal power $Z^{(n)}(a_{n},z_{0};z)$ corresponds to a solution of
(\ref{RiccatiMain}):%
\[
Q_{n}(z)=\frac{\partial_{z}\operatorname*{Re}Z^{(n)}(a_{n},z_{0}%
;z)}{\operatorname*{Re}Z^{(n)}(a_{n},z_{0};z)}%
\]
or $\operatorname*{Re}Z^{(n)}(a_{n},z_{0};z)=e^{A[Q_{n}](z)}$ from where we
obtain (\ref{QviaQ_n}).
\end{proof}

In the next statement we give a generalization of Picard's theorem in the form
(\ref{weyrgen}).

\begin{theorem}
(Picard's theorem) Let $Q_{k}$, $k=1,2,3,4$ be four solutions of
(\ref{RiccatiMain}). Then
\[
\frac{\partial_{\overline{z}}(Q_{1}-Q_{2})+2i\operatorname*{Im}(\overline
{Q}_{1}Q_{2})}{Q_{1}-Q_{2}}+\frac{\partial_{\overline{z}}(Q_{3}-Q_{4}%
)+2i\operatorname*{Im}(\overline{Q}_{3}Q_{4})}{Q_{3}-Q_{4}}%
\]%
\[
-\frac{\partial_{\overline{z}}(Q_{1}-Q_{4})+2i\operatorname*{Im}(\overline
{Q}_{1}Q_{4})}{Q_{1}-Q_{4}}-\frac{\partial_{\overline{z}}(Q_{3}-Q_{2}%
)+2i\operatorname*{Im}(\overline{Q}_{3}Q_{2})}{Q_{3}-Q_{2}}=0.
\]

\end{theorem}

\begin{proof}
Obviously,
\[
(\overline{Q}_{1}+\overline{Q}_{4})+(\overline{Q}_{3}+\overline{Q}%
_{2})-(\overline{Q}_{1}+\overline{Q}_{2})-(\overline{Q}_{3}+\overline{Q}%
_{4})=0.
\]
Multiplying each parenthesis by $1=(Q_{i}-Q_{j})/(Q_{i}-Q_{j})$ we obtain the
equality
\[
\frac{(\overline{Q}_{1}+\overline{Q}_{4})(Q_{1}-Q_{4})}{(Q_{1}-Q_{4})}%
+\frac{(\overline{Q}_{3}+\overline{Q}_{2})(Q_{3}-Q_{2})}{(Q_{3}-Q_{2})}%
\]%
\[
-\frac{(\overline{Q}_{1}+\overline{Q}_{2})(Q_{1}-Q_{2})}{(Q_{1}-Q_{2})}%
-\frac{(\overline{Q}_{3}+\overline{Q}_{4})(Q_{3}-Q_{4})}{(Q_{3}-Q_{4})}=0.
\]
Using
\[
(\overline{Q}_{i}+\overline{Q}_{j})(Q_{i}-Q_{j})=\partial_{\overline{z}}%
(Q_{j}-Q_{i})-\overline{Q}_{i}Q_{j}+Q_{i}\overline{Q}_{j}%
\]
the result is obtained.
\end{proof}

\section{Cauchy's integral theorem}

\begin{theorem}
\label{Th_Cauchy_Int_theorem_Riccati}(Cauchy's integral theorem for the
complex Riccati equation) Let $Q_{0}$ and $Q_{1}$ be bounded solutions of
(\ref{RiccatiMain}) in $\Omega$. Then for every closed curve $\Gamma$ lying in
a simply connected subdomain of $\Omega$,%
\begin{equation}
\operatorname*{Re}\int_{\Gamma}\left(  Q_{1}-Q_{0}\right)  e^{A\left[
Q_{1}-Q_{0}\right]  }dz+i\operatorname*{Im}\int_{\Gamma}\left(  Q_{1}%
-Q_{0}\right)  e^{A\left[  Q_{1}+Q_{0}\right]  }dz=0. \label{CauchyintRiccati}%
\end{equation}

\end{theorem}

\begin{proof}
From theorem \ref{ThFromRiccatiToSchrod} we have that $f=e^{A\left[
Q_{0}\right]  }$ and $u=e^{A\left[  Q_{1}\right]  }$ are solutions of
(\ref{Schrod2}). Now, applying theorem \ref{ThCauchySchr} we obtain%
\[
\operatorname*{Re}\int_{\Gamma}\partial_{z}(\frac{u}{f})dz+i\operatorname*{Im}%
\int_{\Gamma}f^{2}\partial_{z}(\frac{u}{f})dz=0
\]
for every closed curve $\Gamma$ situated in a simply connected subdomain of
$\Omega,$ which gives us the equality%
\[
\operatorname*{Re}\int_{\Gamma}\partial_{z}(e^{A\left[  Q_{1}-Q_{0}\right]
})dz+i\operatorname*{Im}\int_{\Gamma}e^{2A\left[  Q_{0}\right]  }\partial
_{z}(e^{A\left[  Q_{1}-Q_{0}\right]  })dz=0.
\]
From this we obtain the result.
\end{proof}

As a particular case let us analyze the Riccati equation (\ref{RiccatiMain})
with $\nu\equiv0$ which is related to the Laplace equation. If in
(\ref{CauchyintRiccati}) we assume that $Q_{0}\equiv0$, then
(\ref{CauchyintRiccati}) takes the following form%
\[
\int_{\Gamma}Q_{1}e^{A\left[  Q_{1}\right]  }dz=0.
\]
This is obviously valid because if $Q_{1}$ is another bounded solution of the
Riccati equation with $\nu\equiv0$, then according to
\ref{ThFromRiccatiToSchrod} we have that $u=e^{A\left[  Q_{1}\right]  }$ is a
harmonic function and the last formula turns into the equality%
\begin{equation}
\int_{\Gamma}u_{z}dz=0 \label{CauchyintLaplace}%
\end{equation}
($u_{z}$ is analytic).

Now, if in (\ref{CauchyintRiccati}) we assume that $Q_{1}\equiv0$, then
(\ref{CauchyintRiccati}) takes the form%
\[
\operatorname*{Re}\int_{\Gamma}Q_{0}e^{-A\left[  Q_{0}\right]  }%
dz+i\operatorname*{Im}\int_{\Gamma}Q_{0}e^{A\left[  Q_{0}\right]  }dz=0.
\]
Rewriting this equality in terms of the harmonic function $f=e^{A\left[
Q_{0}\right]  }$ we obtain%
\[
\operatorname*{Re}\int_{\Gamma}\frac{f_{z}dz}{f^{2}}+i\operatorname*{Im}%
\int_{\Gamma}f_{z}dz=0
\]
which taking into account (\ref{CauchyintLaplace}) becomes the equality%
\[
\operatorname*{Re}\int_{\Gamma}\frac{f_{z}dz}{f^{2}}=0
\]
or which is the same,%
\begin{equation}
\operatorname*{Re}\int_{\Gamma}\partial_{z}\left(  \frac{1}{f}\right)  dz=0.
\label{vsp23}%
\end{equation}
This equality is a simple corollary of a complex version of the Green-Gauss
theorem (see, e.g., \cite[sect. 3.2]{TutschkeVasudeva}) according to which we
have%
\[
\frac{1}{2i}\int_{\Gamma}\partial_{z}\left(  \frac{1}{f}\right)
dz=\int_{\Omega}\partial_{\overline{z}}\partial_{z}\left(  \frac{1}{f}\right)
dxdy.
\]
For $f$ real the right-hand side is real-valued and we obtain (\ref{vsp23}).

\section{Conclusions}

We have shown that the stationary Schr\"{o}dinger equation in a
two-dimensional case is related to a complex differential Riccati equation
which possesses many interesting properties similar to its one-dimensional
prototype. Besides the generalizations of the famous Euler theorems we have
obtained the generalization of Picard's theorem and the Cauchy integral
theorem for solutions of the complex Riccati equation. The theory of
pseudoanalytic functions has been intensively used.

\begin{acknowledgement}
The authors wish to express their gratitude to CONACYT for supporting this
work via the research project 50424.
\end{acknowledgement}

\end{document}